\documentclass{amsart}
\usepackage{amssymb,amscd}
\usepackage{graphicx}

\usepackage{latexsym}

\usepackage{amsthm}
\newtheorem{theorem}{Theorem}

\newtheorem{defn}{Definition}

\newcommand{\bc}{\mathbb{C}}

\newcommand{\bp}{\mathbb{P}}
\newcommand{\bq}{\mathbb{Q}}
\newcommand{\bz}{\mathbb{Z}}
\newcommand{\modm}{\mathcal{M}}
\newcommand{\modmgn}{\overline{\mathcal{M}}_{g,n}}

\begin{document}

\title{Stable reduction and topological invariants of complex polynomials}
\author{Paul Norbury}
\address{Department of Mathematics and Statistics\\
University of Melbourne\\Australia 3010}
\email{pnorbury@ms.unimelb.edu.au}
\thanks{I would like to thank Brandeis University for its hospitality while this paper was written}

\date{\today}

\begin{abstract}

\noindent  A topological invariant of a polynomial map $p:X\to B$ from a complex surface containing a curve $C\subset X$ to a one-dimensional base is given by a rational second homology class in the compactification of the moduli space of genus $g$ curves with $n$ labeled points $\modmgn$.  Here the generic fibre of $p$ has genus $g$ and intersects $C$ in $n$ points.  In this paper we give an efficient method to calculate this homology class.  We apply this to any polynomial in two complex variables $p :\bc^2\to\bc$ where the $n$ points on a fibre are its points at infinity.
\end{abstract}

\maketitle

\section{Introduction}

A polynomial map $p:X\to B$ comes equipped with a rich topological structure.  The fibres $p^{-1}(b)$ form a family of complex analytic curves parametrised by points in the base $b\in B$.  The generic fibre is a smooth connected topological surface of genus $g$.  A finite set $S=\{b_1,...,b_m\}\subset\bc$ of values in the base have exceptional fibres $p^{-1}(b_i)$ with Euler characteristic strictly greater than the generic Euler characteristic.  On the complement of $S$ the map $p$ defines a surface fibration with structure group the mapping class group of the generic fibre, also known as its geometric monodromy group.  The action of this monodromy group on the homology of the generic fibre is a well-studied topological invariant.  The object of study in this paper is concerned with the geometric monodromy left over after passing to finite covers. For example, if the homological monodromy is trivial, or finite so that a cover is trivial, there can be interesting geometric monodromy remaining.  

An effective way to study the geometric monodromy group uses
the moduli space $\modm_g$ of genus $g$ curves.  This is rationally the classifying space of the mapping class group so it classifies bundles on the complement of the exceptional set.  A polynomial map $p:X\to B$ gives a map from the complement of its exceptional set $B-S\to\modm_g$ by sending $b\in B-S$ to its fibre $p^{-1}(b)$ a smooth genus $g$ curve.  Finite covers of fibrations arise when one tries to extend over the exceptional set to get a map $B\to\overline{\mathcal{M}}_g$ to the compactification which is the moduli space of {\em stable curves} of arithmetic genus $g$.  A stable curve satisfies $g>1$, it has at worst isolated nodal singularities, and any genus zero irreducible component contains at least three singular points \cite{HMoMod}.  Typically, there will be exceptional fibres of $p$ which are not stable.

The main technical issue of how to extend the map 
\[ B-S\to\modm_g\hookrightarrow\overline{\mathcal{M}}_g\]
across the set $S$ of exceptional values is dealt with in Section~\ref{sec:stab} using the process of {\em stable reduction} of a family of curves.  The stable reduction around a singular curve $F_0$ in a family $F_t$ canonically produces a stable curve to replace $F_0$ and induces a canonical topological decomposition along circles of the generic curve in the family.  The stable curve corresponds to pinching along these circles.  The stable reduction is obtained by pulling back the (singular) fibration $p:X\to B$ over a new base $B'\to B$ so that the fibres are reasonably behaved.  In this case reasonably behaved means that the local monodromy around exceptional values is given by Dehn twists around disjoint non-parallel curves, which is equivalent to the homological monodromy being unipotent.  This condition on the monodromy is equivalent to the fibres having at worst isolated nodal singularities.  If we also blow down any genus zero component of a fibre containing fewer than three singular points, in which case the fibre becomes a stable curve, then we have produced a map from $B'\to\overline{\mathcal{M}}_g$, and a corresponding rational homology class $[B']/d\subset H_2(\overline{\mathcal{M}}_g)$ where $d$ is the degree of the cover $B'\to B$.  It is clear that the homology class depends only upon the topological type of the stable reduction and that it is independent of the finite cover.

The process of stable reduction around an exceptional fibre $F_0$ in a family $F_t$ is topologically equivalent to a canonical topological decomposition of the generic curve in the family along circles, together with the degree of the cover.  This process, and hence the decomposition, depends only on a neighbourhood of $F_0$ in the family.  The neighbourhood is encoded in its 3-manifold link, well-defined for arbitrarily small $\epsilon$,
\[ Y=\left\{\bigcup F_t\ |\ |t|=\epsilon\right\}\]
which fibres over $S^1$.  The 3-manifold $Y$ possesses a canonical decomposition known as the JSJ decomposition \cite{JShSei,JohHom}.  The JSJ decomposition canonically decomposes $Y$ along two-dimensional tori.  The intersection of the tori with a  fibre $F_t$ of $Y$ induces a canonical decomposition of the fibre along circles.  
\begin{theorem}  \label{th:decomps}
Given a polynomial map $p:X\to B$, in a neighbourhood of an exceptional fibre the two decompositions of the generic fibre along circles, arising from stable reduction and the JSJ decomposition, coincide.  
\end{theorem}
We have stated the theorem in its simplest possible form although more is true.  The JSJ decomposition also detects further information contained in the stable reduction - the finite cover required during the base change, and the singularities at nodal points of the stable curve.  The singularities at nodal points were only mentioned implicitly in the statement that we should blow down some genus zero components of a fibre.  See Section~\ref{sec:stab} for precise details.  The theorem can also be extended to the situation of a map $p:X\to B$ in the presence of a curve $C\subset X$ that intersects all fibres.   The stronger version of the theorem with proof appears as Theorem~\ref{th:decomps2} in Section~\ref{sec:stab}.  This generalisation is required in our applications and is explained below.  The theorem is simply an interpretation of the methods of Eisenbud and Neumann in \cite{ENeThr} which applies the JSJ decomposition of the 3-manifold boundary neighbourhood of an exceptional fibre to piece together the geometric monodromy.  This ends up being perfectly suited to stable reduction.  The main purpose of the theorem is not only to give a topological construction of the rational homology class associated to $p:X\to B$, but also an efficient way to calculate.  In Section~\ref{sec:blowups} we show the long process involved in stable reduction and in Section~\ref{sec:splice} we show how Theorem~\ref{th:decomps} makes the calculations much more efficient.\\

A curve $C\subset X$ marks points on each fibre of $p:X\to B$.  The generic fibre must intersect $C$ transversally in $n$ points, and in particular the set of exceptional values may enlarge to $S'=\{b_1,...,b_m,b_{m+1},...,b_{m'}\}\subset B$.  As in the case when $C$ is empty, we can characterise an exceptional fibre $p^{-1}(b_i)$ by the property that
the Euler characteristic $\chi(p^{-1}(b_i)-C)$ is strictly greater than the generic Euler characteristic $\chi(p^{-1}(b)-C)$.  The curve $C$ decomposes into irreducible components $C=C_1\cup C_2\cup...\cup C_k$, and we may assume that no $C_i$ lies in a fibre of $p$, since we can just remove such components, so $k\leq n$.  If $k=n$ then the $C_i$ are sections of the map $p:X\to B$ and the generic fibre $p^{-1}(b)$ is an element of $\modm_{g,n}$, the moduli space of genus $g$ curves with $n$ labeled points.  The geometric monodromy fixes the labeled points and $\modm_{g,n}$ is rationally the classifying space of the subgroup of the mapping class group that fixes the $n$ labeled points.  When $k<n$ there is not a well-defined labeling of the $n$ points of intersection on a generic fibre.  In other words, the monodromy acts non-trivially on the intersection of $C$ with the generic fibre, and in particular its action on homology is not unipotent.  There is a finite cover $B'\to B$ with monodromy acting trivially on the intersection of the pull-back of $C$ with the generic fibre, and hence the map $p':X'\to B'$ possesses $n$ sections.  Thus stable reduction, which ensures unipotent monodromy around exceptional fibres, is also used to enable a labeling of the $n$ points of intersection.  As before we get a homology class $[B']/d\subset H_2(\modmgn)$ where $d$ is the degree of the cover $B'\to B$.  The compactification $\modmgn$ consists of stable curves of arithmetic genus $g$ and $n$ labeled points.  Generalising the case $n=0$ to $n>0$, a stable curve has at worst isolated nodal singularities distinct from the labeled points, and any genus zero irreducible component contains at least three points made up of labeled points and singular points.  There is no restriction on genus.
Theorem~\ref{th:decomps} generalises to this situation in which the three-manifold now contains a link $L=Y\cap C\subset Y$.

Since the main purpose of these ideas is to do be able to do calculations, this paper contains examples.  A huge source of examples is supplied by polynomials in two variables $p:\bc^2\to\bc$.  In this case the $n$ points on a generic fibre are its points at infinity.  The stable reduction of a family is obtained by pulling back $p:\bc^2\to\bc$ over a new base $B'\to\bc$ so that the homological monodromy around exceptional fibres and around infinity is unipotent.
This fits in with the picture above if we compactify $\bc^2$ to $X$ and extend $p$ to a map from $X$ to $\bp^1$.  Note that we cannot simply choose $X=\bp^2$ since $p$ does not extend there.  Section~\ref{sec:stab} gives details.  Thus, a polynomial $p :\bc^2\to\bc$ gives a (multi-valued) map from the complement of its exceptional set in its base $\bc-S$ to $\modm_{g,n}$ and stable reduction allows us to extend over the exceptional set and at the same time choose a labeling of the $n$ points at infinity to get a rational homology class in $H_2(\modmgn)$.

The JSJ decomposition has previously been applied to polynomials in two variables in \cite{NeuCom} to introduce a powerful topological invariant of a polynomial $p :\bc^2\to\bc$, its {\em link at infinity} which is a link in the 3-sphere given by the intersection of a generic fibre of $p$ with the boundary of a very large ball in $\bc^2$.  This is analogous to the link of a complex plane curve singularity obtained by intersecting the curve with the boundary of a very small ball around the singularity.  The link at infinity yields topological information about the polynomial such as the genus of the generic curve in the family, the monodromy around exceptional curves in the family, and possible exceptional curves that can arise in the family.  Part of this paper addresses the question of how much the link at infinity of a polynomial can say about its second rational homology class in $\modmgn$. 

An effective description of the homology class associated to $p:X\to B$ uses rational numbers obtained by evaluating the homology class on rational cohomology classes 
\[\delta,\kappa_1,\lambda_1,\psi_i\in H^2(\modmgn,\bq),\ i=1,...,n.\]  
These cohomology classes are defined in Section~\ref{sec:coh} where explicit calculations of the rational numbers are given.

\section{Stable reduction}   \label{sec:stab}

The technique of stable reduction here is taken from \cite{HMoMod}.  
As described in the introduction, generic fibres of a polynomial map $p:X\to B$ with $C\subset X$ are smooth genus $g$ curves containing $n$ distinct points.  Although exceptional fibres are not necessarily stable, the existence of the compactification $\modmgn$ means that as fibres approach an exceptional fibre they converge to a stable curve, so we may replace the exceptional fibre with that stable curve.  This process, known as stable reduction, involves passing to a cover $B'\to B$ in which the local homological monodromy around each exceptional fibre is unipotent.  The pull-back map $p':X'\to B'$ will have the two important properties that every fibre is a stable curve, and that $C'\subset X'$ is the image of $n$ sections $B'\to X'$.  Such a cover exists by the local monodromy theorem which asserts that there exists a power of the local homological monodromy around an exceptional fibre which is unipotent with exponent 2, i.e. $(h^m-1)^2=0$.  There are many such covers since any power of a unipotent element remains unipotent, nevertheless the way in which we define the rational second homology class in the moduli space of stable curves is independent of the cover.  We decompose the cover into cyclic covers of prime degree since they are most easily understood.  The process is best demonstrated via explicit calculations so this is what we will do next.  The reader can consult \cite{HMoMod} for further details.   

\subsection{Branched covers}  \label{sec:blowups}
 
Given a polynomial $p:\bc^2\to\bc$ compactify to a morphism $\hat{p}:X\to\bp^1$.  This is achieved by first compactifying to $\bar{p}:\bp^2\to\bp^1$.  The map $\bar{p}$ is not a morphism (unless it was trivially a coordinate to begin with) or in other words there are points in $\bp^2$ on which $\bar{p}$ takes on every value and hence is ill-defined.  Thus we repeatedly blow up $\bp^2=\bc^2\cup\bp^1_{\infty}$ until the transform $\hat{p}$ of $p$ is well-defined, and in particular there exist rational curves $C_1,...,C_m\subset X-\bc^2$ so that $\hat{p}:C_i\to\bp^1$ is onto, sometimes referred to as horizontal curves or quasi-sections.  If the degree of $\hat{p}:C_i\to\bp^1$ is 1 for each $i$ then $m=n$ and the map $\hat{p}:X\to\bp^1$ has $n$ sections and the generic fibres become points of $\modm_{g,n}$.  If the degree of $\hat{p}:C_i\to\bp^1$ is $d_i>1$ for some $i$ then the inverse map $\bp^1\to C_i$ is a multiply defined $1:d_i$ map.  We must take a finite cover consisting of different choices of labeling the points at infinity meeting $C_i$ to get $d_i$ sections.  We arrange this at the same time that we do stable reduction around a non-stable fibre.  Note that the fibre $\hat{p}^{-1}(\infty)$ is never stable because it has components of multiplicity greater than 1 so we need to replace it.

\begin{figure}[ht] 
	\centerline{\includegraphics[height=2cm]{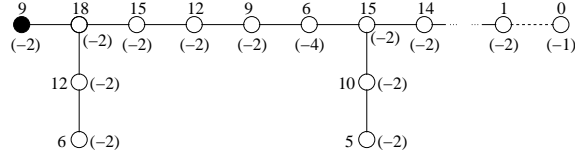}}
	\caption{Fibre of $\hat{p}$ over infinity}
	 \label{fig:res}
\end{figure}

Let $p=(x^2-y^3)^2+xy$.  The generic fibre has genus 4 and 1 point at infinity.   Extend to $\bar{p}:\bp^2\to\bp^1$ and resolve the indeterminate points of $\bar{p}$ (explicitly, change coordinates to $\bar{p}=(z-y^3)^2+yz^5$ and blow up multiple times, or do a single toric blow-up $(y,z)\mapsto(yz,y^2z^3)$, followed by a change of coordinates $(y,z)\mapsto(1-y,z)$ and another toric blow-up $(y,z)\mapsto(y^4z^7,yz^2)$ followed by another change of coordinates and multiple blow-up) to get the curve over infinity represented in Figure~\ref{fig:res} by its dual graph.

In Figure~\ref{fig:res} the nodes denote irreducible components, which are rational curves unless a genus is specified, the edges denote intersections between components, the self-intersection numbers are given in brackets, and the other numbers denote the multiplicity of $\hat{p}$ on each component.  The multiplicities and intersection numbers are related by the fact that a component has zero intersection with the fibre at infinity.  The zero multiplicity on the right-most component denotes the fact that it maps surjectively onto $\bp^1$, and in this case the map has degree 1 so it is a section.  The dark node corresponds to the proper transform of the line at infinity in $\bp^2$.

The Euler characteristic of the generic fibre can be calculated from the fibre over infinity by removing the points of intersection of irreducible components and taking the Euler characteristic.  Thus a node of multiplicity $m$ and valency $v$ is replaced by $m$ copies of a $v$-punctured 2-sphere.  The valency 2 genus 0 components contribute 0 to this sum which is $6+4+5-12-10=-7$ so the generic fibre is a genus 4 curve minus a point at infinity.

Each fibre meets the section exactly once.  The fibre over infinity is not stable since most of its components have multiplicity greater than 1.  All other fibres are stable in this example because the exceptional fibres, which are precisely the singular fibres, contain only nodal singularities. We will demonstrate the process of stable reduction which will replace the fibre at infinity by a stable curve. 

\begin{figure}[ht]  
	\centerline{\includegraphics[height=4.5cm]{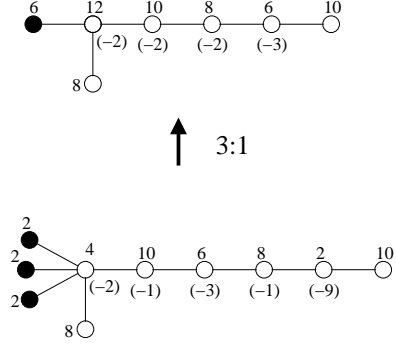}}
	\caption{Example of a 3-fold cover branched over a singular curve}
	\label{fig:3fold}
\end{figure}

Break the branched cover up into a sequence of covers of prime degree.  Take a 3 to 1 branched cover of the base totally ramified at infinity.  The fibre over infinity in the pull-back fibration is not smooth since the original fibre over infinity is not smooth.  For example, the cover $z^3=xy$ is singular above the singularity $(0,0)$, and we can resolve the cover using two intersecting $(-2)$ curves.  More generally we will need to consider $z^p=x^ay^b$ for a $p$-fold cover near the intersection of curves of multiplicities $a$ and $b$.  It turns out that the normalisation of the pull-back fibration is the same as ramifying along curves over infinity of multiplicity not divisible by 3.  For example, in Figure~\ref{fig:3fold} ramify on the curves of multiplicities  8, 10, 8 and 10.  

The rightmost curve of multiplicity 6 in Figure~\ref{fig:3fold} pulls back to a 3-fold cover branched at its two intersection points with the multiplicity 8 and 10 curves.  Its multiplicity is divided by 3.  For the intersecting multiplicity 10 and 8 curves the 3-fold cover ramified along these two is singular.  It is a Hirzebruch-Jung singularity which can be resolved.  In this case it can be resolved by blowing up between the two curves downstairs to get a multiplicity 18 curve which is divisible by 3, so we are reduced to the previous case of a 3-fold cover branched over 2 points, leaving a multiplicity 6 curve.  The multiplicity 12 curve is branched over its two intersection points with multiplicity 8 and 10 and hence 3 copies of the multiplicity 6 curve live in the cover, sharing multiplicity 2 each.

\begin{figure}[ht]  
	\centerline{\includegraphics[height=16cm]{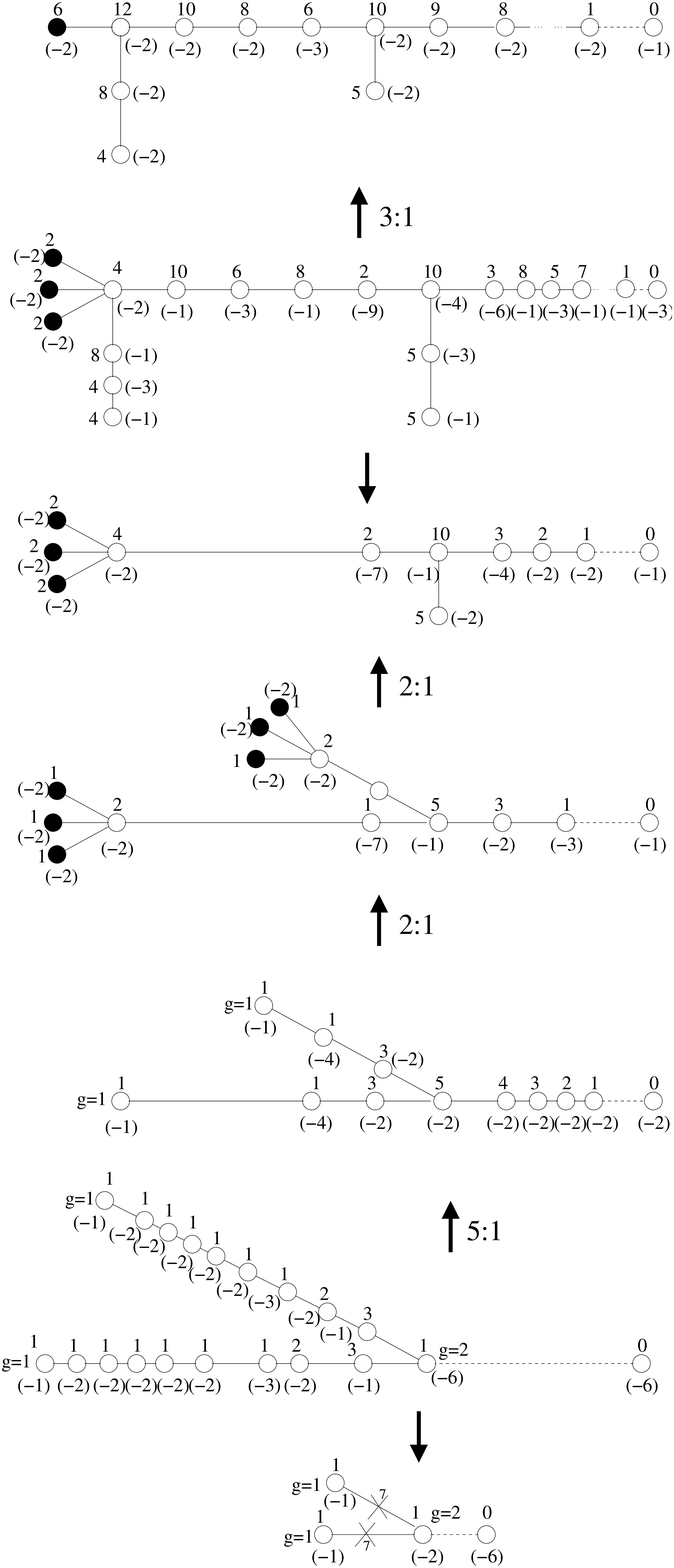}}
	\caption{The fibre over infinity of a 60:1 branched cover is a genus 2 curve meeting two genus 1 curves.}
	\label{fig:stablered}
\end{figure}

Figure~\ref{fig:stablered} shows the fibre over infinity of a sequence of base changes of prime degree combining to give a branched cover $\bp^1\stackrel{60:1}{\to}\bp^1$ that pulls back the fibration so that the fibre over infinity is a stable curve.  Often the resolution of a Hirzebruch-Jung singularity contains a string of curves that can be blown down.  We have not shown the full calculation.  Only once in Figure~\ref{fig:stablered} have we shown the blow-downs preferring to combine steps and hide many intermediate rational curves.

\begin{figure}[ht] 
	\centerline{\includegraphics[height=4cm]{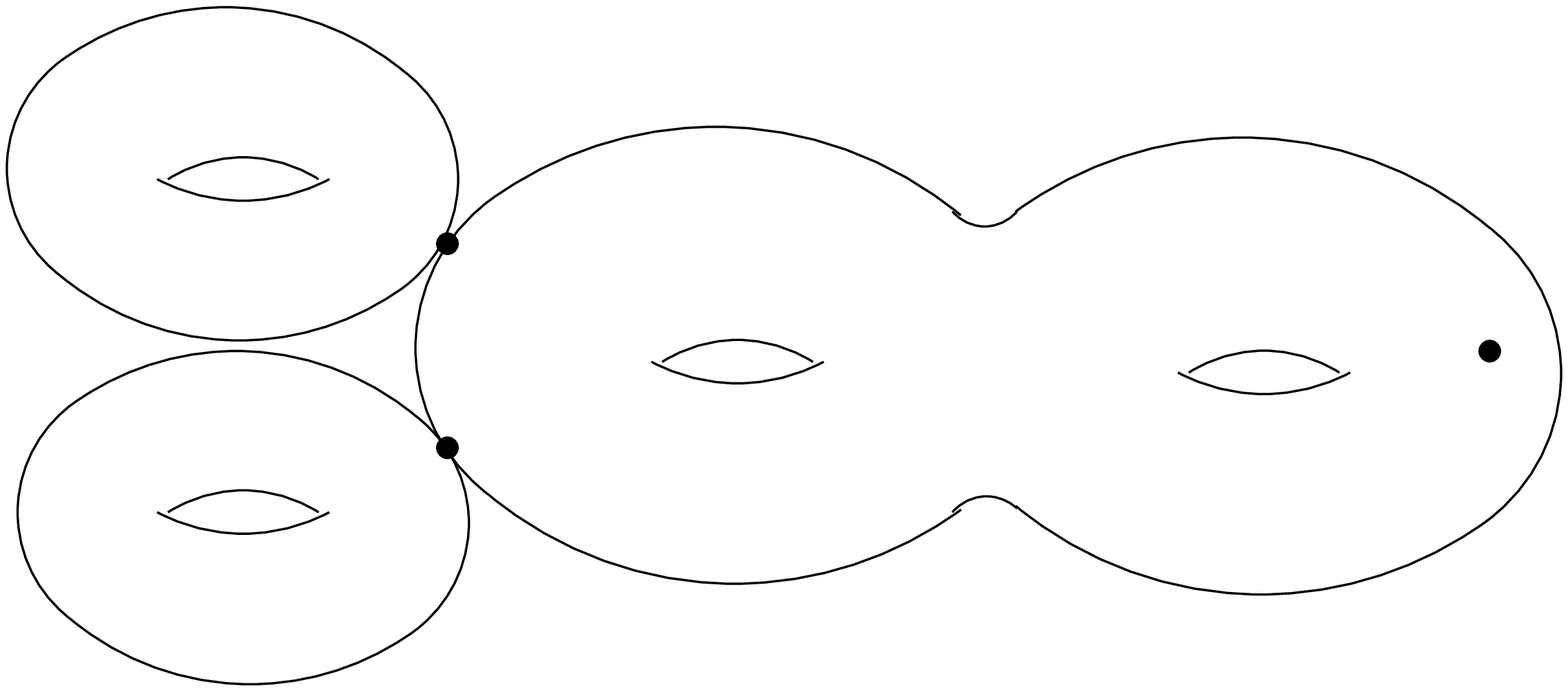}}
	\caption{Stable fibre over infinity.}
	 \label{fig:stableinf}
\end{figure}

The process ends with the stable curve of arithmetic genus 4 given by two genus 1 curves attached to a genus 2 curve as in Figure~\ref{fig:stableinf}.  This stable curve is canonically associated to the family in that any other cover leading to stable reduction will give the same stable curve over infinity.  In the final diagram of Figure~\ref{fig:stablered} the cross with a 7 on an edge means that the intersection of the genus 1 and genus 2 curves is a $\bz_7$ quotient singularity obtained by blowing down a string of six $(-2)$ curves.  In the family of curves, parametrised by $t$, say, with $t=0$ corresponding to infinity, the nodal singularity is given by $xy=t^7$.  Often we stop at the previous step containing strings of $(-2)$ curves, and hence semi-stable fibres, to avoid singularities upstairs.  Of course if we take a further branched cover of degree $d$ say, the nodal singularity would be $xy=t^{7d}$.  It is clear that the factor of 7 is essential and the extra factor of $d$ is the result of taking a less efficient cover.  The essential number is really $7/60$, an invariant of the process.  Similarly, it is $1/60$th of the stable curve in Figure~\ref{fig:stableinf} that strictly replaces the fibre over infinity.  This reflects the fact that the moduli space is an orbifold.

In the example, the curve $C\subset X$ is already a section.  More generally, when there are irreducible components of $C$ that are not sections we need to arrange to have sections in the cover during the stable reduction process.  In each step the $p$-fold cover is totally ramified at infinity and we can freely specify all other ramification.  We choose some of the points of ramification away from infinity to coincide with the points of ramification of the map $C\to\bp^1$ in such a way to remove the ramification.  Note that when the base is $\bp^1$ as for polynomials of two variables, we can also choose the ramification away from infinity so that the cover of the base is again $\bp^1$.

The result of stable reduction produces stable curves to replace non-stable exceptional fibres.  When these new stable curves are included into the family there are usually singularities at nodes of the stable curves.    The degree of the cover of the stable reduction must be remembered in order to get invariant quantities such as the homology class.  The example demonstrates an algorithm for stable reduction that works in general.  Next we will prove Theorem~\ref{th:decomps} and apply it to deduce the stable reduction more efficiently.

\subsection{Canonical decompositions}   \label{sec:splice}
Theorem~\ref{th:decomps} is a special case of a more general theorem which we will now describe.  Consider a polynomial map $p:X\to B$ and a curve $C\subset X$.  For any exceptional fibre $F_0$ define the link of $F_0$ by the pair $(Y,L)=(\partial N(F_0),C\cap\partial N(F_0))$ consisting of a 3-manifold $Y$ containing a 1-manifold $L$, where $N(F_0)$ is a small enough neighbourhood of $F_0$.  The JSJ decomposition decomposes $Y-L$ along finitely many tori.  The manifold $Y-L$ fibres over the circle via $p$, with fibre the complement of $n$ distinct points on the generic fibre of $p$.   Equipped with the more detailed description of stable reduction above, we are ready to restate Theorem~\ref{th:decomps} in its full generality. 

\begin{theorem}  \label{th:decomps2}
The stable reduction of an exceptional fibre $F_0$ of a polynomial map $p:X\to B$ with a curve $C\subset X$ is obtained as follows.  Denote the generic fibre of $p$ by $F$ and let $(Y,L)$ be the link of $F_0$.  Then $Y-L$ fibres over the circle with fibre $F-\{n\ {\rm distinct\ points}\}$ and monodromy $h$ satisfying:

(i) The JSJ decomposition of $Y-L$ along tori induces a decomposition \begin{equation} \label{eq:decomp}
 F-\{n\ {\rm points}\}=F^{(1)}\cup...\cup F^{(k)}
\end{equation} 
of each fibre along circles. 

(ii) The monodromy $h$ preserves the decomposition (\ref{eq:decomp}).  It has finite order on each piece $F^{(i)}$ and acts by a fractional twist along the decomposing circles.

(iii) The stable reduction of $F_0$ in the family uses the cover on which the lift of $h$ is trivial on each $F^{(i)}$.  In particular, the degree of the cover is the lowest common multiple of the finite orders of $h$ on each $F^{(i)}$.

(iv) The quotient singularities at nodes of the stable fibre are detected by the twists of $h$ along the decomposing circles.
\end{theorem}

{\em Remark.}  Consider the case of a fibre $F_0$ that is a nodal curve with isolated singularities and no (-1) rational curves meeting fewer than three other components.  Equivalently the multiplicity of each irreducible component is 1 and each rational curve contains at least two singular points.  More generally, there may be $n$ local sections around $F_0$, and we blow up if necessary to separate  intersection points of the sections on $F_0$.  Then $F_0$ is semi-stable and the stable reduction process is unnecessary.  The nearby generic fibre decomposes $F-\{n\ {\rm points}\}=F^{(1)}\cup...\cup F^{(k)}$ along circles known as vanishing cycles, which correspond to the singular points when pinched.  The monodromy is given by one full positive twist around each vanishing cycle, so in particular it is trivial on each $F^{(i)}$.  A string of $k$ rational curves, each containing exactly two singular points and hence of self-intersection (-2), has $k+1$ vanishing cycles, so the monodromy twists exactly $k+1$ times along the common homology cycle of each vanishing cycle.  If we blow down the string of (-2) rational curves the fibre becomes a stable curve, and the node is a $\bz_{k+1}$ quotient singularity in the ambient surface.   The theorem states that the JSJ decomposition of $Y-L$ is the union of $S^1\times F^{(i)}$ along the boundary tori, and the gluing of the boundary tori encodes the quotient singularities. The tori correspond to intersections of neighbourhoods of irreducible components of $F_0$.  The intersection of two components looks locally like $xy=0$ and $Y$ is locally $\{|x|=\epsilon\}\cup\{|y|=\epsilon\}$ which intersects along the torus $\{|x|=|y|=\epsilon\}$.  
\begin{proof}
Properties (i) and (ii) are proven in \cite{ENeThr}.  There it is shown that $Y-L$ decomposes along tori into Seifert-fibred pieces corresponding to the nodes of valency $>2$ and of positive genus.  The decomposition coincides with the JSJ decomposition so in particular it is canonical.    The fibration $p:Y-L\to S^1$ equips $Y-L$ with an element of $H^{1}(Y-L,\bz)$.  An element of $H^{1}(Y-L,\bz)$ is sometimes known as a {\em multilink} structure on $L$ since it assigns a multiplicity to each component by evaluating along a meridian circle.  Any fibre $F$ of the map $Y-L\to S^1$ has boundary the multilink $L$.  Thus when the multiplicities of all components of $L$ are 1, $F$ is a Seifert surface for $L$, and if a multiplicity of a component $L_i$ is $d_i>1$ then the boundary of $F$ near $L_i$ gives $d_i$ copies of $L_i$.  If $d_i=0$ then one can attach a disk to the $i$th boundary component of $F$ which meets $L_i$ transversally.  The fibre $F$ decomposes as in (\ref{eq:decomp}) and each $F^{(i)}$ is the fibre of a fibration of the corresponding Seifert-fibred piece over the circle.  
The monodromy decomposes and since the monodromy of any Seifert-fibred piece is finite this gives (ii).

Using the definition of a stable curve in terms of its homological monodromy being unipotent, (iii) follows from (i) and (ii), since on a cover where the lift of $h$ is trivial on each $F^{(i)}$, only twisting along separating circles remains.  In the cover, the twisting is integral, not merely fractional.  A precise description of how to calculate the fractional twisting of the monodromy along circles is given in \cite{ENeThr}.  We will show this in calculations below and thus verify (iv).
\end{proof}

The homology class of a polynomial map is trivial if and only if all generic fibres are isomorphic.    Since the stable reduction of such a family is necessarily a trivial fibre bundle, a necessary condition for the generic fibres to be isomorphic is that the stable reduction of any exceptional the fibre is smooth.  Theorem~\ref{th:decomps2} shows that this occurs precisely when the dual graph of the exceptional fibre has only one node of valency $>2$ since this is equivalent to the monodromy having finite order.  Applying this idea to the fibre over infinity of a polynomial $p:\bc^2\to\bc$ gives a severe restriction on polynomials with isomorphic general fibres.  Kaliman classified all polynomials $p:\bc^2\to\bc$ with isomorphic generic fibres \cite{KalPol}.  It seems likely that the classification can be reproven beginning from the restriction just described. 

\subsection{Efficient dual graphs}  \label{sec:effdual}
Many calculations involving the dual graph require only some of the nodes.  Examples are:

(i) the Euler characteristic of the generic fibre;

(ii) a description of branched covers of the dual graph;

(iii) determinants of intersection forms of subgraphs;

(iv) the canonical class supported on the dual graph.\\

\begin{figure}[ht] 
	\centerline{\includegraphics[height=2cm]{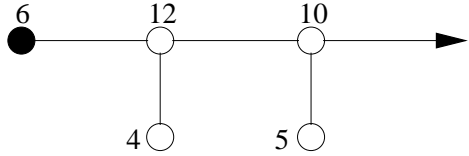}}
	\caption{The efficient dual graph hides valency $2$ genus 0 nodes.}
	 \label{fig:splice}
\end{figure}

An efficient dual graph hides valency 2 rational components, retaining components of valency $\neq 2$ or genus $>0$.  An efficient dual graph for the fibre over infinity of $p=(x^2-y^3)^2+xy$ is given in Figure~\ref{fig:splice}.

We have labeled the multiplicity of the polynomial $p$ on each node.  
The calculation of the Euler characteristic of the generic fibre using the fibre at infinity gets a contribution of zero from the rational valency 2 nodes and hence one needs only the efficient dual graph.  In Figure~\ref{fig:splice} a node of multiplicity $m$ and valency $v$ is replaced by $m$ copies of a $v$-punctured 2-sphere and the Euler characteristic is $6+4+5-12-10=-7$ so the generic fibre is a genus 4 curve minus a point at infinity.\\

A cyclic branched cover along a singular fibre behaves quite generally as in the example in Section~\ref{sec:blowups}.  A cyclic branched cover along a linear chain of genus zero components gives another linear chain of genus zero components.
The example in Section~\ref{sec:blowups} shows this explicitly.
The most interesting behaviour occurs on components intersecting at least two other components, corresponding to nodes of valency $>2$ in the dual graph, and on positive genus components.

By simply understanding the behaviour of a cyclic cover on the nodes of valency $>2$ and of positive genus, we can much more efficiently take cyclic covers of splice diagrams, hiding all chains arising from Hirzebruch-Jung singularities.  The long calculation shown in Section~\ref{sec:blowups} is substantially shortened. \\

\begin{figure}[ht]
	\centerline{\includegraphics[height=2cm]{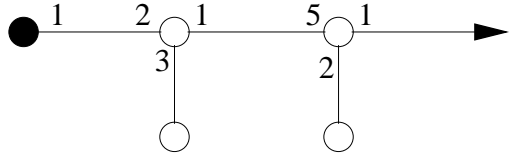}}
	\caption{Splice diagram}
	\label{fig:splice1}
\end{figure}

An efficient dual graph that hides valency $2$ rational nodes, and equipped with enough information to reproduce the full dual graph, is known as a {\em splice diagram.}  Integer weights assigned to the graph supply the information required to reproduce the full dual graph and in particular the multiplicities of each node.  
Determinants of intersection matrices of subgraphs have the property that they are unchanged by blowing up the subgraph.  In particular, to each edge of the efficient dual graph there is a well-defined edge determinant given by the determinant of the intersection matrix of nodes in the full dual graph along that edge, excluding the endpoints.  The efficient dual graph is called a splice diagram when we use a set of weights on half edges.  The weight on a half edge is the determinant of the intersection matrix of the branch of the dual graph disconnected from the half edge if we cut along that edge.  Figure~\ref{fig:splice1} shows Figure~\ref{fig:splice} with weights on half-edges replacing multiplicities on nodes.  The weights give back the whole dual graph. In particular, to get the multiplicity of a node take a path from the node to an arrow and multiply all weights adjacent to, but not on, the path.  The black node has weights 3 and 2 adjacent to the path to the arrow so its multiplicity is $3\times 2=6$.  The closest node to the arrow has weights 5 and 2 adjacent to its path so its multiplicity is $5\times 2=10$.  When there are several arrows do this for each path joining a node to an arrow and add.  

We will not go further into the calculation of the weights, nor will we emphasise the important fact that the splice diagram encodes a 3-manifold containing a link spliced together from Seifert-fibred pieces.  See \cite{ENeThr} for full details.   In this particular case the multiplicities in Figure~\ref{fig:splice} are determined by, and determine, the weights in Figure~\ref{fig:splice1}.\\

Given a polynomial map $p:X\to B$, we may assume the canonical class $K_X$ is supported on fibres of $p$ and sections.  This is because any irreducible component of a divisor $\Sigma\subset X$ is either contained in a fibre of $p$ or it maps onto $B$ and in the latter case we can pass to a cover $B'\to B$ so that the pull-back of $\Sigma$ is a collection of sections, introducing new support only in fibres. 

Let $\cup_iD_i$ denote the support of $K$ and put $K=\sum_ik_iD_i$.  Define $D=\sum_id_iD_i$ where $d_i=-k_i-1$.  Although $D$ resembles a divisor it is not really well-defined since we can add zero times any curve to $K$ and this will change $D$. Nevertheless, $D$ is a useful tool to study the support of $K$.  

The important property of $D$ is that for each $D_i$
\[ D\cdot D_i=\chi(\Sigma)-{\rm valency\ }D_i\]
where the valency of $D_i$ is its valency in the dual graph of $\cup_iD_i$.  This expression vanishes on rational curves of valency 2, precisely those curves hidden in the efficient dual graph.  In particular $K\cdot D$ uses only the multiplicities on the nodes of the efficient dual graph.

For a polynomial $p:\bc^2\to\bc$ the multiplicities of $D$ are easily calculated from the splice diagram.  
The node corresponding to the line at infinity in $\bp^2$, indicated by the black node in the diagrams, has multiplicity 2 since $K$ has multiplicity -3 there.  From any other node $v$ take a path toward the black node to get to the previous node $v'$ of valency $>2$.  Then the multiplicity of $D$  is given recursively by $m(v)=w\cdot m(v')+{\rm edge\ determinant}$ of the edge $E$ joining $v$ and $v'$, and $w$ is the product of the edge weights adjacent to $v$ and not on $E$.   See \cite{NNoOre} for a thorough description of the calculation of $D$ from weights in the splice diagram.  In Figure~\ref{fig:splice1} the multiplicities of $D$ on the two valency 3 nodes are 5 and 3 respectively.  

Another feature of $D$ is that it behaves well under branched covers.  Locally along a curve of ramification $x=0$ the canonical divisor is given by $x^mdx\wedge dy$, and hence a $d$-fold cover $z^d=x$ transforms $z^{md+d-1}dz\wedge dy$, and along unramified curves the multiplicity remains the same.  Thus $D$ is given locally by $x^{-m-1}dx\wedge dy$ transforming to $z^{d(-m-1)}dz\wedge dy$ so the multiplicity simply multiplies by $d$.  

\subsection{Decomposing dual graphs.}

\begin{figure}[ht]  
	\centerline{\includegraphics[height=3cm]{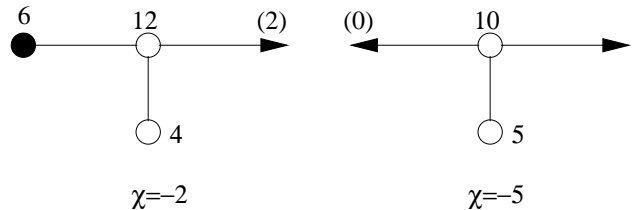}}
	\caption{JSJ decomposition}
	\label{fig:splice3}
\end{figure}

Figure~\ref{fig:splice3} shows how to decompose Figure~\ref{fig:splice} into pieces and how the Euler characteristic decomposes.  Each node of valency $>2$ corresponds to a Seifert-fibred piece of the JSJ decomposition of $Y-L$.  The $(2)$ on the arrow indicates that a fibre of the left piece consists of two identical disjoint components.  The greatest common divisor of the $(0)$ on an arrow and the default value of $(1)$ on the right arrow is 1 which is the number of components of the fibre of the right piece.  The numbers $(2)$ and $(0)$ are calculated using the weights of the splice diagram in the same way as multiplicities are calculated, by taking a path from the edge to an arrow and multiply all weights adjacent to, but not on, the path.  Thus the left piece has fibre two genus 1 surfaces each with one boundary component and the right piece has fibre a genus 2 surface with three boundary components.  The genus 4 generic fibre with one point removed decomposes along two circles as shown in Figure~\ref{fig:stable}.
By pinching along these circles, or vanishing cycles, one gets the stable curve in Figure~\ref{fig:stableinf}.

\begin{figure}[ht] 
	\centerline{\includegraphics[height=4cm]{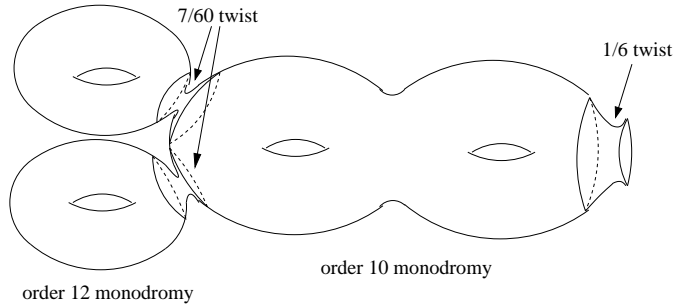}}
	\caption{The monodromy acts with finite order on each piece and fractionally twists along the decomposing annuli.}
	 \label{fig:stable}
\end{figure}

This demonstrates (i) of Theorem~\ref{th:decomps2}.  The decomposition in Figure~\ref{fig:splice3} gives the JSJ decomposiiton of $Y-L$ into pieces $Y^{(i)}$ each fibring over the circle and containing a multilink $L^{(i)}$.  The numbers $(2)$ and $(0)$ in Figure~\ref{fig:splice3} are the multiplicities of the multilink components.  The monodromy on the fibre of $Y-L$ decomposes to give the monodromy of the fibre of $Y^{(i)}-L^{(i)}$.  The latter has finite order given by the multiplicity of the valency $>2$ component.  This uses the fact that the monodromy on Seifert-fibred spaces is completely understood.  Thus the multiplicity is 12 on the left component - this consists of an order 6 map on the once-punctured torus and a swapping of the two components - and the multiplicity is 10 on the right component.  A 60:1 cover will pull back the monodromy to be trivial on each component and thus this will give the stable reduction of the fibre over infinity as seen in the calculation in Section~\ref{sec:blowups}.  Thus (iii) and part of (ii) Theorem~\ref{th:decomps2} are seen.  

It remains to understand the twisting of the monodromy along the annuli in Figure~\ref{fig:stable}.  Upstairs, in the 60:1 cover we expect some whole number of Dehn twists along the annuli, so downstairs we expect fractional twists.  This calculation is subtler and it is completely solved in terms of the JSJ decomposition as Theorem 13.1 in \cite{ENeThr}.  It is subtle only because we are choosing to emphasise the more intuitive multipicities of the polynomial $\hat{p}$ on each component of the efficient dual graph, rather than the weights of the splice diagram.  The fractional twisting requires the weights of the splice diagram.  It is given along an annulus corresponding to the edge $E$ by
\[{\rm twist}(h|A)=\frac{-d_E}{l\cdot l'}\Delta_E\]
where $d_E$ is the number of annuli corresponding to $E$,
$l$ and $l'$ are the multiplicities of the neighbouring nodes, and $\Delta_E$ is the determinant of the intersection matrix of the nodes from the full dual graph lying on the edge. In the example above, $d_E=2$, coresponding to the 2 annuli, $l\cdot l'=12\cdot 10=120$ and $\Delta_E=\left|\begin{array}{ccc}-2&1&0\\1&-2&1\\0&1&-2\end{array}\right|=-7=5\cdot 1-2\cdot 3\cdot 2\cdot 1$.  (We have given the calculation of the edge determinant twice, the second using the splice diagram weights in Figure~\ref{fig:splice1}, not explained here further.)  Thus the fractional twisting is $7/60$ and upstairs in the 60:1 cover it is a Dehn twist of order 7.  This tells us that locally the node of the stable curve is given by $xy=t^7$ and there is a $\bz_7$ quotient singularity in the ambient surface.

\section{Cohomology classes}  \label{sec:coh}
In this section we evaluate the homology class associated to $p:X\to B$ on rational cohomology classes 
\[\delta,\kappa_1,\lambda_1,\psi_i\in H^2(\modmgn,\bq),\ i=1,...,n.\]  
Since we consider 1-dimensional bases $B$, we will abuse notation and treat the cohomology classes as numbers $\delta, \kappa_1,\lambda_1,\psi_i\in\bq$. 

A non-trivial homology class in $H_2(\modmgn)$ is detected by the non-vanishing of at least one of the cohomology classes above.  This follows from work of Wolpert \cite{WolHom} where he shows that $\kappa_1+\sum_{i=1}^n\psi_i$ is a positive multiple of the K\"ahler class $\omega$.  The integral of $\omega$ over any non-trivial homology class is positive and hence at least one of $\kappa_1$ or $\psi_i$ is non-zero on the homology class.

There are obstructions to a dual graph appearing as the fibre at infinity
of a polynomial $p:\bc^2\to\bc$.  The obstructions are not complete so there are dual graphs where it is not known if they appear as the fibre at infinity of a polynomial $p:\bc^2\to\bc$.  We do not know if the positivity condition $\kappa_1+\sum_{i=1}^n\psi_i>0$ gives a new obstruction on dual graphs.

\subsection{Definitions}
The class $\delta\in H^2(\modmgn)$ is represented by the boundary divisor $\modmgn-\modm_{g,n}$.  Thus it counts the number of singularities of a polynomial $p:X\to B$.  The contribution from a stable fibre is its number of nodes, where a node locally given by $xy=t^k$ contributes $k$ to the count.  Equivalently, if we insist that the total space $X$ be non-singular and thus allow semi-stable fibres then the contribution from a semi-stable fibre is simply its number of nodes.  The contribution from an unstable fibre is calculated by first passing to the stable reduction.  In particular, the contribution is not necessarily an integer.  
\begin{defn}
Let $p':X'\to B'$ be the semi-stable reduction of $p:X\to B$ equipped with $n$ sections, for $X'$ a smooth surface that $d$-fold covers $X$.  Define
\begin{eqnarray*}\delta&=&\frac{\#{\rm\ singularities\ }p'}{d}\\
\kappa_1&=&2\delta+\frac{3\sigma(X')}{d}\\
\lambda_1&=&\frac{\delta}{4}+\frac{\sigma(X')}{4d}\\
\psi_i&=&-\frac{s_i(B')\cdot s_i(B')}{d},\quad i=1,...,n
\end{eqnarray*}
where $\sigma(X')$ is the signature of $X'$ and $s_i(B')\cdot s_i(B')$ is the self-intersection of the $i$th section of $X'\to B'$.
\end{defn}
It is not hard to see that if we take a further cover the number of singularities multiplies by the degree of the cover and hence $\delta$ is well-defined.  It is subtler that the signature behaves the same way.
The signature of a branched cover is in general not easily calculated.  The definition above is given for simplicity whereas for the purpose of calculations we give a better definition below in which the classes $\kappa_1$ and $\lambda_1$ in families $\kappa_m$ and $\lambda_m$ for $m\in\bz^+$.   The relation between the two definitions of $\lambda_1$ is due to Smith \cite{SmiLef} and the others easily follow.

The forgetful map ${\overline{\mathcal{M}}_{g,n+1}}\to\modmgn$ is defined by forgetting the $(n+1)$st point and possibly blowing down rational components of some fibres.  It defines a universal bundle over $\modmgn$ equipped with $n$ sections $s_i:\modmgn\to{\overline{\mathcal{M}}_{g,n+1}}$.  Over ${\overline{\mathcal{M}}_{g,n+1}}$ define  the vertical canonical bundle to be the complex line bundle 
\[ \gamma=K_{{\overline{\mathcal{M}}_{g,n+1}}}\otimes\pi^*K_{\modmgn}^{-1}\to{\overline{\mathcal{M}}_{g,n+1}}.\]
The bundle $\gamma$ is used to define the classes $\kappa_1,\lambda_1,\psi_i$ as follows.

Define the Hodge classes
\[\lambda_m=c_m(\pi_*(\gamma))\in H^{2m}(\modmgn).\]
The push-forward sheaf $\pi_*(\gamma)$ is a rank $g$ vector bundle over $\modmgn$ best understood in terms of the fibres which are the $g$-dimensional vector spaces of holomorphic 1-forms on the curve associated to a point  of $\modmgn$.  Over a stable curve, one uses 1-forms holomorphic outside the singular set that have at worst simple poles at the singular points with residues summing to zero.

Define the Mumford-Morita-Miller classes
\[\kappa_m=\pi_![c_1(\gamma)^{m+1}]\in H^{2m}(\modmgn).\] 
where $\pi_!:H^k({\overline{\mathcal{M}}_{g,n}})\to H^{k-2}(\modmgn)$ is the umkehr map, or Gysin homomorphism, obtained by integrating along the fibres. 

For each $i=1,...,n$ pull back the line bundle $\gamma$ to $s_i^*\gamma=\gamma_i\to\modmgn$ and define 
\[ \psi_i=c_1(\gamma_i)\in H^2(\modmgn).\]

The classes are related \cite{HMoMod} by
\[ \lambda_1=\frac{\kappa_1+\delta}{12}.\]

\subsection{Calculations}
The main purpose of this paper is to enable calculations of the homology class of a polynomial map $p:X\to B$ with 1-dimensional base $B$ and in particular, calculations of the rational numbers
$\delta, \kappa_1,\lambda_1,\psi_i$.  We will calculate these rational numbers for infinite families of polynomial maps $p:\bc^2\to\bc$.  We first do the calculations for the specific example $p=(x^2-y^3)^2+xy$.

The contribution from the fibre over infinity to $\delta$ comes from the 
two nodes of the stable fibre over infinity, each counted with multiplicity 7.   The stable reduction is obtained from a 60:1 cover hence the contribution to $\delta$ is $\frac{7}{30}$.  This is the fractional twisting of the geometric monodromy and we have given it a cohomological interpretation.  There are 8 finite singularities so 
\[\delta=8\frac{7}{30}.\]

\begin{figure}[ht] 
	\centerline{\includegraphics[height=3cm]{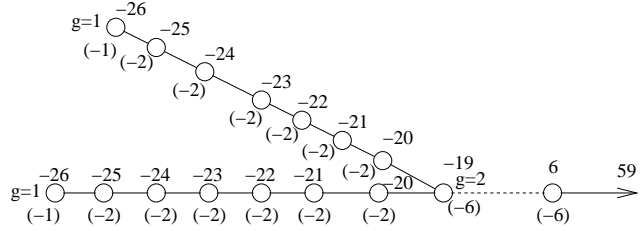}}
	\caption{Canonical divisor}
	\label{fig:canonical}
\end{figure}

To calculate $\kappa_1$ we need to understand the canonical class on the surface $X'$ which gives the stable family $X'\to B'$.  The stable fibre over infinity is given in Figure~\ref{fig:canonical}
where the numbers in brackets are self-intersection numbers as usual and the multiplicities are now the multiplicities of the canonical class. 

We deduce the canonical class by the fact that it can be supported on the infinite fibre plus 59 copies of the fibre of $p$.  This is because the canonical class can be supported on the divisor at infinity $X-\bc^2$ since it can be supported on $\bp^2-\bc^2$ and blowing up introduces support only on the exceptional divisiors.  A branched cover introduces canonical class along the ramification set so the 60:1 cover means 59 copies of the fibre of $p$ contributes to the canonical class and the rest remains supported on the fibre over infinity.  Apply the adjunction formula $K\cdot\Sigma=-\chi(\Sigma)-\Sigma\cdot\Sigma$ to a fibre to see $K\cdot F=-6$ so the section is given multiplicity 6 in the canonical class. To calculate the other multiplicities in Figure~\ref{fig:canonical}, similarly apply the adjunction formula to irreducible components of the fibre over infinity. 

Thus 
\[60\kappa_1=K_{X'/B'}^2=(K_{X'}+2F)^2=274.\]
and
\[\lambda_1=\frac{\kappa_1+\delta}{12}=\frac{4\frac{17}{30}+8\frac{7}{30}}{12}=\frac{16}{15}.\]
Note that $\lambda_1$ times the degree of the cover is an integer, so the relation $\lambda_1=(\kappa_1+\delta)/12$ puts a mod 12 condition on $\kappa_1+\delta$.
Finally $\psi_1$ is given by $K_{X'/B'}\cdot s_1[B']=-B'\cdot B'$ by the adjunction formula.  Thus
\[\psi_1=\frac{6}{60}=\frac{1}{10}.\]
The self-intersection $(-6)$ used to calculate $\psi_1$ is encoded in the fractional twisting of the geometric monodromy.  It gives the fractional twisting of the generic fibre on the annulus around the marked point.  In terms of the efficient dual graph, or splice diagram, the annulus $A$ corresponds to the edge $E$ joining the fibre over infinity to the section.  Thus $\psi_1$ is given by
\[\psi_1=-{\rm twist}(h|A)=\frac{d_E}{l}=\frac {1}{10}\]
where $d_E=1$ is the number of annuli corresponding to $E$,
and $l$ is the multiplicities of the node of $E$.  The formula for the twist is Theorem 13.5 in \cite{ENeThr}.

\subsubsection{Calculations in ${\overline{\mathcal{M}}_{g,1}}$.}
The example $p=(x^2-y^3)^2+xy$ is a particular case of the general class of polynomials with one point at infinity.  The splice diagram for the general polynomial with one point at infinity is given in Figure~\ref{fig:stable}.
\begin{figure}[ht] 
	\centerline{\includegraphics[height=2cm]{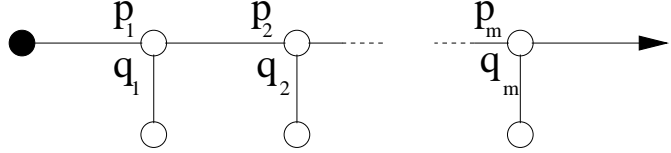}}
	\caption{Family in ${\overline{\mathcal{M}}_{g,1}}$.}
	 \label{fig:stable}
\end{figure}

The stable reduction of the fibre over infinity is easily calculated from Theorem~\ref{th:decomps}.   We refer to the node with weights $p_i,q_i$ as the $i$th node and the edge between the $i$th and $(i+1)$th nodes as the $i$th edge.  The $i$th node corresponds to $\prod_{j>i}q_j$ components of Euler characteristic $p_i+q_i-p_iq_i$, each meeting $q_i$ of the $\prod_{j>i-1}q_j$ components corresponding to the $(i-1)$th node.  The monodromy permutes the $\prod_{j>i}q_j$ components and acts internally with order $p_iq_i$, so its total order on the components corresponding to the $i$th node is $p_iq_i\prod_{j>i}q_j$ which is the multiplicity of the polynomial on that node.  Thus the stable reduction uses a cover of order $\prod p_iq_i$ (or more efficiently $gcd(p_1,...,p_n)\prod q_i$.)  

From the fractional twisting along annuli we see that the quotient singularity at each annulus above the $i$th edge has order $p_iq_iq_{i+1}-p_{i+1}$, and the self-intersection of the section in the stable reduction is $-\prod_{i<n}p_iq_i$.  We calculate inequalities for the rational numbers $\delta, \kappa_1,\lambda_1$ since their maximum value is taken when the fibre of $p$ over infinity is the only unstable curve, and the rational numbers $\delta, \kappa_1,\lambda_1$ are upper-semicontinuous.  

We have no good interpretation of the messy expressions for $\delta$ and $\kappa_1$ so we show them only for the cases $m=1$ and $m=2$.  For $m=1$ the Euler characteristic of the generic fibre is
\[\chi=p_1+q_1-p_1q_1,\quad(m=1).\]
The stable reduction of the fibre over infinity uses a degree $p_1q_1$ cover and results in a smooth fibre over infinity. 
To calculate $\delta$ we use a result of Suzuki \cite{SuzPro}
\begin{equation}  \label{eq:suz} 
\sum_{c\in\bc}(\chi_c-\chi)=1-\chi
\end{equation}
where $\chi_c=\chi(p^{-1}(c))$ which expresses the Euler characteristic of the generic fibre as a sum of the number of vanishing cycles near exceptional fibres.  If all fibres over finite values $c\in\bc$ are stable, which is the generic case that gives an upper bound on $\delta$, then $\chi_c-\chi$ counts the number of nodes in the fibre $p^{-1}(c)$.  Thus
\[\delta\leq 1-\chi,\quad(m=1).\]
The fractional twisting along the boundary annulus is given by $-1/p_mq_m$ so
\[\psi_1=\frac{1}{p_mq_m}.\]
If all finite fibres are stable, the canonical class upstairs in the $p_1q_1$-fold cover is $K=(-2-\chi)F+(-1-\chi)C$ where $F$ is the generic fibre and $C$ is the section.  The fractional twisting $-1/p_1q_1$ shows that $C\cdot C=-1$ in the cover.  Thus $K\cdot K=-(1+\chi)^2+2(1+\chi)(2+\chi)=(1+\chi)(3+\chi)$ and
\[\kappa_1\leq \frac{K\cdot K-4(\chi+1)}{p_1q_1}=\frac{\chi^2-1}{p_1q_1},\quad(m=1).\]
When $m=2$, the stable reduction of the fibre over infinity uses a cover of order $p_1q_1p_2q_2$ and consists of an irreducible component with Euler characteristic $1+p_2+q_2-p_2q_2$ plus $q_2$ irreducible components with Euler characteristic $1+p_1+q_1-p_1q_1$, joined along $q_2$ nodes with quotient singularities of order $p_1q_1q_2-p_2$.  Now
\[\chi= q_1 q_2 + p_1 q_2 + p_2 - p_1 q_1 q_2 - p_2 q_2,\quad(m=2)\]
and
\[\delta\leq 1-\chi-\frac{1}{p_1q_1}+\frac{q_2}{p_2},\quad(m=2)\]
uses (\ref{eq:suz}) to get an upper bound for nodes away from infinity, plus $1/(p_1q_1p_2q_2)$ times the $q_2$ nodes over infinity, each counted with multiplicity $p_1q_1q_2-p_2$.

To calculate $\kappa_1$ we need a more efficient method than has been used so far.  We use the divisor $D$ introduced in Section~\ref{sec:effdual}.  The expression
\begin{equation}  \label{eq:kdotk}
 K\cdot K=-K\cdot D-\sum_i(\chi(C_i)+C_i\cdot C_i)
\end{equation}
consists of $K\cdot D$ which is independent of rational valency 2 nodes in the dual graph, and $\sum_i(\chi(C_i)+C_i\cdot C_i)$ which does depend on rational valency 2 nodes in the dual graph, however it behaves well on the semi-stable reduction of a family.  A semi-stable curve contains strings of $(-2)$ rational curves in place of the quotient singularities at nodes of stable curves.  The expression $\sum_i(\chi(C_i)+C_i\cdot C_i)$ vanishes on $(-2)$ rational curves and hence (\ref{eq:kdotk}) uses only the multiplicities of $K$ on the irreducible components of the stable fibres, and on sections.  From this we calculate
\[\kappa_1\leq 2-2\chi+\frac{1}{p2q2}-p1q1-p2q2+\frac{p1}{q1}+\frac{q1}{p1}+\frac{p2}{q2},\quad(m=2).\]

\subsubsection{Calculations in ${\overline{\mathcal{M}}_{0,n}}$.}

When $p:\bc^2\to\bc$ has rational fibres $\lambda=0=\kappa+\delta$ and in particular we get equalities instead of inequalities for $\kappa$ and $\delta$, i.e. the generic behaviour always occurs, since if the sum of two upper-semicontinuous functions is continuous then the two functions are both continuous.

When the $n$ points at infinity are sections, no branched covering is required for stable reduction.  We use the classification in \cite{NNoRat} to see that either the homology class is trivial or $\kappa_1=1-n$, $\delta=n-1$, $\psi_i=1$, $i=1,...,n-1$ and $\psi_n=n-3$.  The homology class depends only on $n$.  There are many inequivalent polynomials for each $n$ showing that the homology class is a rather course invariant.  The homological monodromy is trivial in both cases when $\kappa_1=0$ or $\kappa_1=1-n$ so the homology class is useful to distinguish these two cases.

Kaliman \cite{KalRat} classified all polynomials with rational fibres and one fibre isomorphic to $\bc^*$.  We find that $\kappa_1=-1$ and $\delta=1$ on the whole family.

\begin{small}

\end{small}

\end{document}